\documentclass[12pt]{article}
\usepackage{graphicx,color}
\usepackage{lscape}
\usepackage{epsfig}
\usepackage{subfigure}
\usepackage{hhline,amssymb,epsfig,amsmath,array,amsthm}
\begin{document}

\newcommand{\ga}{{\alpha}}
\newcommand{\gb}{{\beta}}
\newcommand{\gc}{{\chi}}
\newcommand{\gd}{{\delta}}
\newcommand{\gep}{{\epsilon}}
\newcommand{\gvare}{\varepsilon}
\newcommand{\gga}{{\gamma}}
\newcommand{\gk}{{\kappa}}
\newcommand{\gl}{{\lambda}}
\newcommand{\gm}{{\mu}}
\newcommand{\gn}{{\nu}}
\newcommand{\go}{{\omega}}
\newcommand{\gp}{{\pi}}
\newcommand{\gph}{{\phi}}
\newcommand{\gvp}{{\varphi}}
\newcommand{\gps}{{\psi}}
\newcommand{\gth}{{\theta}}
\newcommand{\gr}{{\rho}}
\newcommand{\gs}{{\sigma}}
\newcommand{\gt}{{\tau}}
\newcommand{\gx}{{\xi}}
\newcommand{\gz}{{\zeta}}
\newcommand{\gE}{{\Upsilon}}
\newcommand{\gGa}{{\Gamma}}
\newcommand{\gL}{{\Lambda}}
\newcommand{\gO}{{\Omega}}
\newcommand{\gP}{{\Pi}}
\newcommand{\gPh}{{\Phi}}
\newcommand{\gPs}{{\Psi}}
\newcommand{\gTh}{{\Theta}}
\newcommand{\gS}{{\Sigma}}
\newcommand{\gX}{{\Xi}}
\newcommand{\gU}{{\Upsilon}}

\newcommand{\fa}{{\bf a}}
\newcommand{\fb}{{\bf b}}
\newcommand{\fc}{{\bf c}}
\newcommand{\fd}{{\bf d}}
\newcommand{\fe}{{\bf e}}
\newcommand{\ff}{{\bf f}}
\newcommand{\fk}{{\bf k}}
\newcommand{\fm}{{\bf m}}
\newcommand{\fn}{{\bf n}}
\newcommand{\fo}{{\bf o}}
\newcommand{\fp}{{\bf p}}
\newcommand{\fq}{{\bf q}}
\newcommand{\fr}{{\bf r}}
\newcommand{\fs}{{\bf s}}
\newcommand{\ft}{{\bf t}}
\newcommand{\fu}{{\bf u}}
\newcommand{\fv}{{\bf v}}
\newcommand{\fw}{{\bf w}}
\newcommand{\fx}{{\bf x}}
\newcommand{\fy}{{\bf y}}
\newcommand{\fz}{{\bf z}}

\newcommand{\fA}{{\bf A}}
\newcommand{\fB}{{\bf B}}
\newcommand{\fC}{{\bf C}}
\newcommand{\fD}{{\bf D}}
\newcommand{\fE}{{\bf E}}
\newcommand{\fF}{{\bf F}}
\newcommand{\fG}{{\bf G}}
\newcommand{\fI}{{\bf I}}
\newcommand{\fJ}{{\bf J}}
\newcommand{\fK}{{\bf K}}
\newcommand{\fL}{{\bf L}}
\newcommand{\fM}{{\bf M}}
\newcommand{\fN}{{\bf N}}
\newcommand{\fO}{{\bf O}}
\newcommand{\fP}{{\bf P}}
\newcommand{\fQ}{{\bf Q}}
\newcommand{\fR}{{\bf R}}
\newcommand{\fS}{{\bf S}}
\newcommand{\fT}{{\bf T}}
\newcommand{\fU}{{\bf U}}
\newcommand{\fV}{{\bf V}}
\newcommand{\fW}{{\bf W}}
\newcommand{\fX}{{\bf X}}
\newcommand{\fY}{{\bf Y}}
\newcommand{\fZ}{{\bf Z}}

\newcommand{\cB}{{\cal B}}
\newcommand{\cE}{{\cal E}}
\newcommand{\cG}{{\cal G}}
\newcommand{\cS}{{\cal S}}
\newcommand{\cF}{{\cal F}}
\newcommand{\cT}{{\cal T}}
\newcommand{\cD}{{\cal D}}
\newcommand{\cI}{{\cal I}}
\newcommand{\cK}{{\cal K}}
\newcommand{\cL}{{\cal L}}
\newcommand{\cM}{{\cal M}}
\newcommand{\cN}{{\cal N}}
\newcommand{\cU}{{\cal U}}
\newcommand{\cW}{{\cal W}}
\newcommand{\cX}{{\cal X}}
\newcommand{\cY}{{\cal Y}}
\newcommand{\cZ}{{\cal Z}}

 \newcommand{\fga}{\mbox{\boldmath $\alpha$}}
  \newcommand{\fgb}{\mbox{\boldmath $\beta$}}
  \newcommand{\fgd}{\mbox{\boldmath $\delta$}}
 \newcommand{\fgg}{\mbox{\boldmath $\gamma$}}
  \newcommand{\fgl}{\mbox{\boldmath $\lambda$}}
  \newcommand{\fgm}{\mbox{\boldmath $\mu$}}
 \newcommand{\fgs}{\mbox{\boldmath $\sigma$}}
 \newcommand{\fgth}{\mbox{\boldmath $\theta$}}
  \newcommand{\fgt}{\mbox{\boldmath $\tau$}}
  \newcommand{\fgz}{\mbox{\boldmath $\zeta$}}
  \newcommand{\fgps}{\mbox{\boldmath $\psi$}}
  \newcommand{\fge}{\mbox{\boldmath $\eta$}}
  \newcommand{\fgga}{\mbox{\boldmath $\gamma$}}
  \newcommand{\fgv}{\mbox{\boldmath $\varphi$}}
  \newcommand{\fgp}{\mbox{\boldmath $\pi$}}
  \newcommand{\fvare}{\mbox{\boldmath $\varepsilon$}}
  \newcommand{\fgD}{{\bf \Delta}}
  \newcommand{\fgG}{{\bf \Gamma}}
  \newcommand{\fgL}{{\bf \Lambda}}
  \newcommand{\fgPh}{{\bf \Phi}}
  \newcommand{\fgPi}{{\bf \Pi}}
  \newcommand{\fgPs}{{\bf \Psi}}
  \newcommand{\fgS}{{\bf \Sigma}}
  \newcommand{\fgTh}{{\bf \Theta}}
  \newcommand{\fgO}{{\bf \Omega}}
  \newcommand{\fgX}{{\bf \Xi}}
  \newcommand{\fgU}{{\bf \Upsilon}}

\newcommand{\Mean}{{\mbox{E}}}
\newcommand{\Cov}{{\mbox{cov}}}
\newcommand{\Var}{{\mbox{var}}}
\newcommand{\Corr}{{\mbox{corr}}}
\newcommand{\diag}{{\mbox{diag}}}
\newcommand{\prob}{{\mbox{Pr}}}

\newcommand{\Nul}{{\bf 0}}
\newcommand{\One}{{\bf 1}}
\newcommand{\Bd}{\B_{\bullet}}
\newcommand{\Xd}{\X_{\bullet}}
\newcommand{\Zd}{\Z_{\bullet}}
\newcommand{\Nor}{N}
\newcommand{\f}{\textbf}
\newcommand{\Real}{{\cal R}}
\newcommand{\Natural}{{\cal N}}
\newtheorem{thm}{Theorem}
\newtheorem{cor}[thm]{Corollary}
\newtheorem{lem}[thm]{Lemma}
\newtheorem{prop}[thm]{Proposition}
\newtheorem{ax}{Axiom}
\theoremstyle{definition}
\newtheorem{defn}{Definition}[section]
\theoremstyle{remark}
\newtheorem{rem}{Remark}[section]
\newtheorem*{notation}{Notation}
\newcommand{\secref}[1]{\S\ref{#1}}
\newcommand{\bysame}{\mbox{\rule{3em}{.4pt}}\,}
\newcommand{\A}{\mathcal{A}}
\newcommand{\B}{\mathcal{B}}
\newcommand{\st}{\sigma}
\newcommand{\XcY}{{(X,Y)}}
\newcommand{\SX}{{S_X}}
\newcommand{\SY}{{S_Y}}
\newcommand{\SXY}{{S_{X,Y}}}
\newcommand{\SXgYy}{{S_{X|Y}(y)}}
\newcommand{\Cw}[1]{{\hat C_#1(X|Y)}}
\newcommand{\G}{{G(X|Y)}}
\newcommand{\PY}{{P_{\mathcal{Y}}}}
\newcommand{\X}{\mathcal{X}}
\newcommand{\wt}{\widetilde}
\newcommand{\wh}{\widehat}

\newcommand{\bY}{\mathbf{Y}}
\newcommand{\bV}{\mathbf{V}}
\newcommand{\bG}{\mathbf{G}}
\newcommand{\bU}{\mathbf{U}}
\newcommand{\bM}{\mathbf{M}}
\newcommand{\bH}{\mathbf{H}}
\newcommand{\bD}{\mathbf{D}}
\newcommand{\bI}{\mathbf{I}}
\newcommand{\by}{\mathbf{y}}
\newcommand{\bx}{\mathbf{x}}
\newcommand{\bv}{\mathbf{v}}
\newcommand{\bz}{\mathbf{z}}
\newcommand{\btheta}{\boldsymbol{\theta}}
\newcommand{\bpsi}{\boldsymbol{\psi}}
\newcommand{\bX}{\mathbf{X}}
\newcommand{\bepsilon}{\boldsymbol{\epsilon}}
\newcommand{\bbeta}{\boldsymbol{\beta}}

\begin{center}
{\LARGE\bf Simultaneous Model Selection and Estimation for Mean and Association Structures with Clustered Binary Data}
\medskip
\medskip
\medskip


B{\scriptsize Y} XIN GAO\\
{\it Department of Mathematics and Statistics, York University, Toronto, Onatrio \\ Canada M3J 1P3} \\
xingao@mathstat.yorku.ca\\
{\scriptsize AND} GRACE Y. YI\\
{\it Department of Statistics and Actuarial Sciences, University of Waterloo, Waterloo, Canada} \\
yyi@uwaterloo.ca\\

\end{center}

\newpage

\centerline{ABSTRACT}

\noindent This paper investigates the property of the penalized
estimating equations when both the mean and association structures 
are modelled. To select variables for the mean and association structures sequentially, we propose a hierarchical penalized generalized estimating equations (HPGEE2) approach. The first set of penalized estimating equations is solved for the selection of significant mean parameters. Conditional on the selected mean model, the second set of penalized estimating equations is solved for the selection of significant association parameters. The hierarchical approach is designed to accommodate possible model constraints relating the inclusion of covariates into the mean and the association models. This two-step penalization strategy enjoys a compelling advantage of easing computational burdens compared to solving the two sets of penalized equations simultaneously.
HPGEE2 with a smoothly clipped
absolute deviation (SCAD) penalty is shown to have the oracle
property for the mean and association models. The asymptotic behavior of the penalized estimator under this hierarchical approach is established. An efficient two-stage penalized weighted least square algorithm is developed to implement the proposed method. The empirical performance of the proposed 
HPGEE2 is demonstrated through Monte-Carlo studies and the analysis of a
clinical data set.


\begin{center}
1. I{\scriptsize NTRODUCTION}
\end{center}
Clustered binary data arise in many application areas of the biological and social sciences. For example, the disease status of members within a family is correlated due to the sharing of common genetic factors and environmental background. The objectives of statistical analysis often focus on modeling of the relationship between the response and explanatory variables and the association among response measurements within clusters. First-order generalized estimating equations (GEE) have been developed for the analysis of longitudinal and other types of clustered data (Liang \& Zeger, 1986). This approach does not require the specification of the joint distribution of the responses and is widely used to make inferences about the regression parameters on the marginal means. Recently, there has been increasing interest in the inferences on the association models. Second-order GEEs have been proposed to model the association between the responses. Prentice (1988) developed second-order GEEs with emphasis on the estimation of correlation parameters. Fitzmaurice \& Laird (1993) proposed a model which parameterizes the association in terms of conditional odds ratios. Lipsitz, Laird \& Harrington (1991), Liang, Zeger \& Qaqish (1992), Carey, Zeger \& Diggle (1993), Molenberghs \& Lesaffre (1994), Lang \& Agresti (1994), and Fitzmaurice \& Lipsitz (1995) proposed models that parameterizes the association in terms of marginal odds ratio. Other discussion can be found in Yi \& Cook (2002), Yi, He \& Liang (2009, 2011) and He \& Yi (2011), among others. 

When we model the mean and the association structures, we may face a large collection of covariates in which some of them are not important to feature the mean and the association structures of the responses. It is desirable to have a model selection approach which can select the significant covariates for both the mean and association models. In the literature, most variable selection procedures have been developed for the selection of covariates that influence the mean responses (e.g.,  Tibshirani, 1996, Fan \& Li, 2002, 2004; Qu \& Li, 2006; Garcia, Ibrahim \& Zhu, 2010; Wang, 2011). In a different direction, many other methods have been proposed to select covariance structures for multivariate Gaussian or binary data under the framework of graphical models (e.g., Meinshausen \& Buhlmann, 2006; Yuan \& Lin, 2007; Friedman, Hastie \&  Tibshirani, 2007). In Bondell, Krishna \& Ghosh (2010) and Ibrahim, et al. (2011), fixed and random effects selection is considered in mixed effects models. So far, however, there has been little development on variable selections on mean and association models simultaneously.

When both the mean and association structures are modeled, we may be interested in model selection for mean model only, or for association model only, or for both mean and association models together.  Suitable constraints may be imposed on the selection procedures to reflect individual analysis objectives.  For example, in some applications, practitioners may add constraints that the exclusion of certain covariates  from the mean model conceptually implies the exclusion of those covariates from the association model. This poses a technical question of how to properly incorporate such model constraints into selection and estimation procedures for the mean and association models. To address all of these concerns, we propose a hierarchical model selection strategy. In principle, we first select the covariates for the mean model, and then conditional on the selected mean model, we proceed to select variables for the association model. However, an important issue arises, as the selection procedure for the mean model involves estimation of the association parameters. To overcome this difficulty, we use the estimators obtained from the usual unpenalized GEE2 equations as the initial estimates. We first apply a one-step penalization on the mean parameters, then based on the penalized mean estimator, apply a one-step penalization on the association parameters. The penalized estimators obtained from the proposed HPGEE2 using a SCAD penalty is shown to have the oracle
property. The asymptotic behavior of the penalized estimator under this hierarchical approach is investigated. An efficient two-stage penalized weighted least square algorithm is developed for the implementation of the proposed method.

The rest of the article is organized as follows. Section 2 and 3 introduce the  development of hierarchical penalized GEE2 method. Section 4 establishes the theoretical properties of the proposed method. Section 5 provides the details of the implementation of the algorithm. In Sections 6 and 7, empirical performance of the proposed penalized
GEE2 is demonstrated through Monte-Carlo studies and the analysis of a
clinical data set.
  
\begin{center}
2. N{\scriptsize OTATION AND
 MODEL FORMULATION}
\end{center}

Suppose $n$ clusters are randomly selected for the study. For cluster
$i=1,\dots,n,$ let $Y_i=(Y_{i1},\dots,Y_{in_i})^T$ be an
$n_i\times 1$ binary response vector with mean $E(Y_i)=\mu_i,$ and let
$\phi_{ijk}$ be the odds ratio between responses $Y_{ij}$ and
$Y_{ik}$ $(1\leq j<k\leq n_i)$ defined by
$$
\phi_{ijk}=\frac{P(Y_{ij}=1,Y_{ik}=1)P(Y_{ij}=0,Y_{ik}=0)}{P(Y_{ij}=0,Y_{ik}=1)P(Y_{ij}=1,Y_{ik}=0)}
.$$
We consider a marginal model as follows: 1) $g(\mu_{ij})=x_{ij}^T\beta,$ where $g(.)$ is a known link
function, $x_{ij}$ is a $p\times 1$ vector of explanatory
variables associated with $Y_{ij},$ and the $\beta$ are regression
coefficients to be estimated; 2) $\log \phi_{ijk}=z_{ijk}^T\alpha,$ where $z_{ijk}$ is a
$q\times 1$ vector of covariates which specifies the form of the
association between $Y_{ij}$ and $Y_{ik},$ and $\alpha$ is a
$q\times 1$ vector of association parameters to be estimated.

To ease computation, in contrast to the method of Prentice (1988), Carey, Zeger \& Diggle (1993) proposed the alternating logistic regression method. The strategy is to estimate $\alpha$ using the $_{n_i}C_2$ conditional events, $Y_{ij}$ given $Y_{ik}$ and the covariates, for $1\leq j <k \leq n_i$ with an appropriate offset:
\begin{align}
\label{offset} \text{logit}
P(Y_{ij}=1|Y_{ik}=y_{ik}, x_{ij}, z_{ijk})=(z_{ijk}^T\alpha) y_{ik}+\log
\bigl(\frac{\mu_{ij}-\nu_{ijk}}{1-\mu_{ij}-\mu_{ik}-\nu_{ijk}}\bigr),
\end{align}
where  $\nu_{ijk}=E(Y_{ij}Y_{ik}|x_{ij},x_{ik},z_{ijk})$ for $1\leq j<k\leq n_i,$ and $\nu_i=(\nu_{ijk}, 1\leq j < k \leq n_i)^T.$   Let $\zeta_{i}$ denote the $_{n_i}C_2$-vector
with elements
$$
\zeta_{ijk}=E(Y_{ij}|Y_{ik}=y_{ik}, x_{ij}, z_{ijk})=\text{logit}^{-1}\{(z_{ijk}^T\alpha)y_{ik}+\log
\bigl(\frac{\mu_{ij}-\nu_{ijk}}{1-\mu_{ij}-\mu_{ik}-\nu_{ijk}}\bigr)\},
$$
and let $R_i$ be the vector of residuals with elements
$
R_{ijk}=Y_{ij}-\zeta_{ijk}.
$
Let $S_i$ denote the diagonal
matrix with diagonal element $\zeta_{ijk}(1-\zeta_{ijk})$, and
$T_i$ denote the matrix $\partial
\zeta_i/\partial \alpha^T.$ We define
$
A_i=Y_i-\mu_i,$
$B_i=\Cov(Y_i|x_i,z_i),$ and 
$C_i=\partial {\mu_i}/{\partial {\beta}^T},
$
where $x_i=(x_{ij},1\leq j \leq n_i)^T$ and $z_i=(z_{ijk},1\leq j < k \leq  n_i)^T.$
The two sets of equations for the mean and association parameters are given by 
$$
U_{\beta}=\sum_{i=1}^n C_i^T B_i^{-1}A_i=0,
$$
and
$$
U_{\alpha}=\sum_{i=1}^n T_i^T S_i^{-1}R_i=0.
$$
The two equations are solved iteratively by alternating between two steps: 1) Given the association parameter $\alpha$, estimate $\beta$ as a parameter in a marginal logistic regression using the first set of equations; 2) For a given $\beta$, estimate the odds ratio parameter $\alpha$ using a logistic regression of $Y_{ij}$ on each $Y_{ik}$ $(k>j)$ with offset that involves $\mu_{ij}$ and $\nu_{ijk}.$ Let $\hat{\beta}_A$ and $\hat{\alpha}_A$ denote the resultant estimators. Under regularity conditions, $\hat{\beta}_A$ and $\hat{\alpha}_A$ are consistent estimators. 

\begin{center}
3. M{\scriptsize ETHODOLOGY}
\end{center}

\begin{center}
3.1 O{\scriptsize BJECTIVES AND ISSUES}
\end{center}
If we consider variable selection and estimation in GEE2 setting, there could be three different scenarios: a) We can select the mean model while the association parameters are assumed to be appropriate and therefore are consistently estimated; b) We are interested in the selection of association model while the mean parameters are assumed to be appropriate and therefore consistently estimated; c) We are interested in the simultaneous selection of mean and association parameters. Therefore, it is desirable to have a unified approach which can accommodate all of these scenarios. Furthermore, in the last scenario, we need to consider the relationship between the two sets of covariates $x_{ij}$ and $z_{ijk},$ which can be classified as (1) $\{x_{ij},x_{ik}\}\cap z_{ijk}=\Phi;$
(2) $\{x_{ij},x_{ik}\}= z_{ijk};$
(3) $\{x_{ij},x_{ik}\}\cap z_{ijk}\neq\Phi.$
Scenario (2) is actually a special case of scenario (3). In the latter two scenarios, the covariates for the mean model are also considered as potential covariates for the association model. Then there may or may not exist constraints relating the selection of covariates. For example, some models may require that if $x_{ij}$ is significant predictor with nonzero $\alpha_x$ in the association model, it implies that $x_{ij}$ is a significant predictor with nonzero $\beta_x$ in the mean model. For instance, in the study of disease occurrence rates in a family, if a genetic marker is not significant for the disease occurrence rates for each individual, scientists may speculate that this marker cannot influence the association between the disease occurrence among the family members. Therefore, for the underlying true model, it would be sensible to add such constraints as $\beta_{x0}=0$ implying $\alpha_{x0}=0.$ This would then create an additional issue for developing procedures  for selecting mean and association models. To address all of these considerations, we propose a hierarchical model selection strategy. We first select the covariates for the mean model and then conditional on the selected mean model, we proceed to select the association model. However, some conceptual issues arise, as the selection procedure for the mean model involves estimation of the association parameters. But such difficulty can be solved if we use the estimators for the unpenalized GEE2 equations as the initial estimates and then apply a one-step penalization, first on the mean parameters, then based on the penalized mean estimator, apply a one-step penalization on the association parameters.

\begin{center}
3.2 M{\scriptsize ODEL SELECTION AND ESTIMATION}
\end{center}

We propose to perform the model selection for mean parameters
by using:
\begin{align}
\label{equation1}
U_{\beta}(\beta,\hat{\alpha}_A)-n  p_{\lambda}'(|\beta|)\odot \text{sign}(\beta)=0,
\end{align}
where $p_{\lambda}'(|\beta|)=(p_{\lambda}'(|\beta_1|),\dots,p_{\lambda}'(|\beta_p|))^T$ is a $p$-dimensional vector of penalty functions, $\text{sign}(\beta)=(\text{sign}(\beta_1),\dots,\text{sign}(\beta_p))^T$ with $\text{sign}(t)=I(t>0)-I(t<0),$ and $I(.)$ is the indicator function. The notation $\odot$ denotes the component-wise product.
Let $\hat{\beta}$ be the penalized estimator for the mean parameter, obtained by solving Equation (\ref{equation1}).

 Next we perform the model selection for the association parameters by using the following function:
\begin{align}
\label{equation2}
U_{\alpha}(\alpha,\hat{\beta})-n  p_{\lambda}'(|\alpha|)\odot \text{sign}(\alpha)=0,
\end{align}
where $p_{\lambda}'(|\alpha|)=(p_{\lambda}'(|\alpha_1|),\dots,p_{\lambda}'(|\alpha_q|))^T$ is a $q$-dimensional vector of penalty functions, and $\text{sign}(\alpha)=(\text{sign}(\alpha_1),\dots,\text{sign}(\alpha_q))^T.$

Remark 1:
As we start with the consistent estimators $\hat{\beta}_A$ and $\hat{\alpha}_A,$ one-step penalization with proper penalty leads to consistent and sparse estimators $\hat{\beta}$ and $\hat{\alpha}$, with many of the elements estimated to be zero. Furthermore, we are able to show that such penalized estimators enjoy the oracle property. That is, with probability tending to one, the procedure selects the correct sub mean model and the correct sub association model. Our sequential two-stage approach allows the selection of covariates for the association to be conditional on the result from the selection of mean model. This property has an advantage to incorporate into selection procedures the constraints concerning the mean and association covariates. For example, if there is a model constraint that $\beta_{x0}=0$ leads to $\alpha_{x0}=0,$ then if $\hat{\beta}_{x}=0,$ the procedure would automatically set $\hat{\alpha}_{x}=0,$ so that the solution satisfies the constraint. Due to the stagewise nature,  this method allows users to perform only model selection on either the mean parameters or the association parameters or both. If both steps utilize nonzero penalties, we jointly selects the mean and association models.  If there are no constraints relating the mean and the association models, an alternative approach could be simultaneously solving the two sets of penalized estimating equations. This alternative procedure would involve alternating between the two penalized equations until convergence, which will be computationally more intensive than the proposed sequential two-stage method.

Remark 2:
In terms of the choice of penalty function, there are many penalty functions available.  As the LASSO penalty, $p_{\lambda}
(|\theta_l|)=\lambda|\theta_{l}|,$ increases
linearly with the size of its argument, it leads to biases for the
estimates of nonzero coefficients. To attenuate such estimation
biases, Fan and Li (2001) proposed the SCAD penalty. The penalty
function satisfies $p_{\lambda}(0)=0,$ and its first-order
derivative is
\begin{align*}
p_{\lambda}'(\theta)=\lambda\{I(\theta\leq \lambda)+
\frac{(a\lambda-\theta)_{+}}{(a-1)\lambda}I(\theta>
\lambda)\},\,\,\text{for}\,\,\theta\geq 0,
\end{align*}
where $a$ is some constant, usually set to $3.7$ (Fan and Li,
2001), and $(t)_{+}=t I(t>0)$ is the hinge loss function.
The SCAD penalty is a quadratic spline function with knots at
$\lambda$ and $a\lambda$. It is singular at the origin which
ensures the sparsity and continuity of the solution. The penalty
function does not penalize as heavily as the $L_1$ penalty
function on parameters with large values. It has been shown that the likelihood estimation with the SCAD penalty
not only selects the correct set of significant covariates, but
also produces parameter estimators as efficient as if we know the
true underlying sub-model (Fan \& Li, 2001). Namely, the estimators have the so-called oracle property. However, it has not been investigated if the oracle
property is also enjoyed by hierarchical penalized estimating equations with the SCAD
penalty. Furthermore, the asymptotic behavior of the hierarchical penalized GEE2 estimators needs to be investigated. In the next section, we will address these issues.

\begin{center}
4. T{\scriptsize HEORY}
\end{center}

Let $\beta_0$ and $\alpha_0$ denote the true value of $\beta$ and $\alpha,$ respectively. Define $U_{\beta\beta}(\beta,\alpha)=\partial U_{\beta}(\beta,\alpha)/\partial \beta^T,$  $U_{\beta\alpha}(\beta,\alpha)$ $=\partial U_{\beta}(\beta,\alpha)/\partial \alpha^T,$  $U_{\alpha\beta}(\beta,\alpha)=\partial U_{\alpha}(\beta,\alpha)/\partial \beta^T,$ and  $U_{\alpha\alpha}(\beta,\alpha)=\partial U_{\alpha}(\beta,\alpha)/\partial \alpha^T.$ Let
$
H_{\beta\beta}(\beta,\alpha)=$$E_{\beta_0,\alpha_0}$ $\{- U_{\beta\beta}(\beta,\alpha)\},
$
$
H_{\beta\alpha}(\beta,\alpha)=$$E_{\beta_0,\alpha_0}\{- U_{\beta\alpha}(\beta,\alpha)\},
$
$
H_{\alpha\beta}(\beta,\alpha)$ $=E_{\beta_0,\alpha_0} \{U_{\alpha\beta}(\beta,\alpha)\},
$
and
$
H_{\alpha\alpha}(\beta,\alpha)$ $=E_{\beta_0,\alpha_0} $$\{- U_{\alpha\alpha}(\beta,\alpha)\}
$
for the case with $n=1.$
We define the matrix of $H(\beta,\alpha)$ as
$$
\left(
    \begin{array}{cc}
      H_{\beta\beta}(\beta,\alpha), & H_{\beta\alpha}(\beta,\alpha) \\
      H_{\alpha\beta}(\beta,\alpha), & H_{\alpha\alpha}(\beta,\alpha) \\
    \end{array}
  \right).
$$
Denote the covariance matrix of the estimating equations for $n=1$ as
$
V_{\beta\beta}(\beta,\alpha)=$$\text{Cov}_{\beta_0,\alpha_0}$ $\{ U_{\beta}(\beta,\alpha)\},
$
$
V_{\beta\alpha}(\beta,\alpha)=V_{\alpha\beta}(\beta,\alpha)^T=$$\text{Cov}_{\beta_0,\alpha_0}\{U_{\beta}(\beta,\alpha),U_{\alpha}(\beta,\alpha))^T\},
$
$
V_{\alpha\alpha}(\beta,\alpha)=$$\text{Cov}_{\beta_0,\alpha_0}$ $\{ U_{\alpha}(\beta,\alpha\}.
$
Furthermore, we define the matrix of $V(\beta,\alpha)$ as
$$
\left(
    \begin{array}{cc}
      V_{\beta\beta}(\beta,\alpha), & V_{\beta\alpha}(\beta,\alpha) \\
      V_{\alpha\beta}(\beta,\alpha), & V_{\alpha\alpha}(\beta,\alpha) \\
    \end{array}
  \right).
$$

For notational convenience, let $s$ denote the set of $j$ such that $\beta_{j0}\neq 0$ and $s^c$ denote the set of $j$ such that $\beta_{j0}= 0.$ Let $v$ denote the set of $j$ such that $\alpha_{j0}\neq 0$ and $v^c$ denote the set of $j$ such that $\alpha_{j0}= 0.$ To establish the asymptotic properties of the proposed penalized estimators, we assume certain regularity conditions which are listed in Appendix A.

\begin{thm}
\label{thm1}
Let $G(\beta,\hat{\alpha}_A)=U_{\beta}(\beta,\hat{\alpha}_A)-n  p_{\lambda}'(|\beta|)\odot \text{sign}(\beta),$ and 
$G_j(\beta,\hat{\alpha}_A))$ be its $j$th element, $j=1,\dots,p.$ Given the SCAD penalty function $p_{\lambda}(\theta),$ if $\lambda_n\rightarrow 0,$ and $\sqrt{n}\lambda_n \rightarrow
\infty$ as $n\rightarrow \infty,$ then there exist a solution $\hat{\beta}$ such that $G(\hat{\beta},\hat{\alpha}_A)=0,$ and $||\hat{\beta}-\beta_0||=O_p(n^{\frac 1 2}).$ Furthermore, we have
$$
\lim_{n\rightarrow \infty} P(\hat{\beta}_{j}=0)=1,
$$ for all $j$ such that $\beta_{j0}=0,$
\end{thm}

Theorem 1 establishes the existence of consistent estimator $\hat{\beta}$
 to $G(\beta,\hat{\alpha}_A)=0.$  Furthermore, the estimator has the property of setting the nonsignificant mean parameter to zero with probability tending to one. Next, we  establish the asymptotic distribution of the penalized mean estimator $\hat{\beta}$.
Let $\beta_s=(\beta_{10},\dots,\beta_{p'0})^T$ be the subset of nonzero mean parameter, $
\Sigma_1=\text{diag}\{p''_{|\lambda_n|}(\beta_{10}),\dots,p''_{\lambda_n}(|\beta_{p'0}|)\},
$
and
$
b_1=(p'_{\lambda_n}(\beta_{10})\text{sign}(\beta_{10}),\dots,$ $p'_{\lambda_n}\text{sign}(\beta_{p'0}))^T.
$
Let
$V_{ss}$ denote the submatrix of $V(\beta_0,\alpha_0)$ corresponding to the index subset $s.$

\begin{thm}
\label{thm2}
Given the SCAD penalty function $p_{\lambda}(\theta),$
if $\lambda_n\rightarrow 0$ and $\sqrt{n}\lambda_n\rightarrow
\infty,$ as $n \rightarrow \infty,$ then the sub-vector of the root-n consistent estimator  $
\hat{\beta}_s$ has the following asymptotic
distribution:
$$
\sqrt{n}(H_{ss}+\Sigma_1)\{\hat{\beta_s}-\beta_{s0}+(H_{ss}+\Sigma_1)^{-1}b_1\}\rightarrow
N\{0,V_{ss}\}, \, \text{as}\, n\rightarrow \infty.
$$
\end{thm}

Now we investigate the properties of the penalized estimator for the association parameter. Let $J(\hat{\beta},\alpha)=U_{\alpha}(\hat{\beta},\alpha)-n  p_{\lambda}'(|\alpha|)\odot \text{sign}(\alpha),$ where
$J_j(\hat{\beta},\alpha)$ denotes its $j$th element, $j=1,\dots,q.$ Let $J_v$ denote the set of penalized equations corresponding to indices set $v,$ and $U_{\alpha_v}$ denote the subset of $U_{\alpha}$ corresponding to indices set $v.$ Denote the inverse of $H(\beta_0,\alpha_0)$ as
$$
\left(
    \begin{array}{cc}
      H^{\beta\beta} & H^{\beta\alpha} \\
      H^{\alpha\beta} & H^{\alpha\alpha} \\
    \end{array}
  \right).
$$
As $H_{\beta\alpha}=0,$ it can be shown that $H^{\beta\alpha}=0,$ $H^{\beta\beta}=H_{\beta\beta}^{-1},$ and $H^{\alpha\alpha}=H_{\alpha\alpha}^{-1}.$ Let $H_{ss}$, $H_{sv},$ $H_{vs},$ and $H_{vv}$ represent the submatrices of $H(\beta_0,\alpha_0),$ let $H^{ss}$, $H^{sv},$ $H^{vs},$ $H^{vv}$ represent the submatrices of  $H(\beta_0,\alpha_0)^{-1},$ and let $V_{ss}$, $V_{sv},$ $V_{vs},$ and $V_{vv}$ represent the submatrices of $V(\beta_0,\alpha_0)$ corresponding to the subset of indices $s$ and $v.$

\begin{thm}
\label{thm3}
Given $\hat{\beta}$ as the solution to the first set of penalized equation $G(\beta,\hat{\alpha}_A)$ with the SCAD penalty function, if $\lambda_n\rightarrow 0,$ and $\sqrt{n}\lambda_n \rightarrow
\infty$ as $n\rightarrow \infty,$ then there exists a solution $\hat{\alpha}$ such that $J(\hat{\beta}, \hat{\alpha})=0,$ and $||\hat{\alpha}-\alpha_0||=O_p(n^{\frac 1 2}).$ Furthermore, we have
$
\lim_{n\rightarrow \infty} P(\hat{\alpha}_{j}=0)=1,
$ for all $j \in v.$
\end{thm}

Next, we need to establish the asymptotic joint distribution of the penalized estimator for the mean and the association parameters. Denote by $\alpha_v=(\alpha_{10},\dots,\alpha_{q'0})^T$ the subset of nonzero association parameters, 
$
\Sigma_2=\text{diag}\{p''_{\lambda_n}(|\alpha_{10}|),\dots,p''_{\lambda_n}(|\alpha_{q'0}|)\},
$
and
$
b_2=(p'_{\lambda_n}(\beta_{10})\text{sign}(\alpha_{10}),\dots,$ $p'_{\lambda_n}\text{sign}(\alpha_{q'0}))^T.
$

\begin{thm}
\label{thm4}
Given the SCAD penalty function, if $\lambda_n\rightarrow 0$ and $\sqrt{n}\lambda_n\rightarrow
\infty,$ as $n \rightarrow \infty,$ then the sub-vectors of the root-n consistent estimators  $
\hat{\beta}_s$ and $
\hat{\alpha}_v$ have the following joint asymptotic
distribution:
$$
\sqrt{n}\left(
          \begin{array}{c}
            \hat{\beta}_s-\beta_{s0} \\
            \hat{\alpha}_v-\alpha_{v0} \\
          \end{array}
        \right)\rightarrow N \left(
                               \begin{array}{cc}
                                 \left(
                                   \begin{array}{c}
                                     -b_1 \\
                                     -b_2 \\
                                   \end{array}
                                 \right),
                                  &B \left(
                                      \begin{array}{cc}
                                        V_{ss} & V_{sv} \\
                                        V_{vs} & V_{vv} \\
                                      \end{array}
                                    \right)B^T
                                   \\
                               \end{array}
                             \right)
, \, \text{as}\, n\rightarrow \infty,
$$ where
$$B= \left(
                                      \begin{array}{cc}
                                        H_{ss}+\Sigma_1 & 0 \\
                                        H_{vs} & H_{vv}+\Sigma_2 \\
                                      \end{array}
                                    \right)^{-1}.$$
\end{thm}

\begin{center}
5. I{\scriptsize MPLEMENTATION: PENALIZED RE-WEIGHTED LEAST SQUARE LGORITHM}
\end{center}

For notational convenience, Let $C=(C_1^T,\dots,C_n^T)^T,$ $A=(A_1^T,\dots,A_n^T)^T,$ and
$B$ be block diagonal with $B_i$
as the diagonal elements.  Let
$T=(T_1^T,\dots,T_n^T)^T,$ $R=(R_1^T,\dots,$ $R_n^T)^T,$ and
$S$ be block diagonal with $S_i$
as the diagonal elements. Let $u=C^{T}B^{-1}A,$ $D=C^TB^{-1}C,$ $u^*=T^{T}S^{-1}R,$ and $D^*=T^TS^{-1}T.$   The initial estimates are $\hat{\beta}^{(0)}=\hat{\beta}_A,$ and $\hat{\alpha}^{(0)}=\hat{\alpha}_A.$ We employ the following two steps for selection and estimation purposes.
\begin{itemize}
\item Step1 : Selection and estimation for the mean model using penalized estimating equation (\ref{equation1}).

\begin{itemize}
\item  Outer loop: Based on the current estimate $\hat{\beta}^{(t)},$ we compute the updated estimate $\hat{\beta}^{(t+1)}.$ We iterate the update, $t=1,2, 3, \dots,$ until $||\hat{\beta}^{(t)}-\hat{\beta}^{(t+1)}||\leq \epsilon,$ the prespecified tolerance level. Using modified Fisher scoring method, we obtain
\begin{align}
\begin{split}
\hat{\beta}^{(t+1)}=&\hat{\beta}^{(t)}-\{\sum_{i=1}^n
C_i(\hat{\beta}^{(t)})^T B_i(\hat{\beta}^{(t)},\hat{\alpha}_A)^{-1}
C_i(\hat{\beta}^{(t)})\}^{-1}\\
&\{\sum_{i=1}^n C_i(\hat{\beta}^{(t)})^T B_i(\hat{\beta}^{(t)},\hat{\alpha}_A)^{-1}
A_i(\hat{\beta}^{(t)})+n\sum_{l=1}^p p'_{\lambda_n}(|\hat{\beta}^{(t)}_l|)\text{sign}(\hat{\beta}^{(t)}_l)\}.
\end{split}
\end{align}
Update the modified dependent variable
$Z(\hat{\beta}^{(t)})=C^T(\hat{\beta}^{(t)})\hat{\beta}^{(t)}-A(\hat{\beta}^{(t)}).$
Then the $\hat{\beta}^{(t+1)}$ solves the penalized weighted linear
regression of responses $Z$ on design matrix $C$ with weight matrix $B^{-1}.$ That
implies $\hat{\beta}^{(t+1)}$ minimizes the following objective
function:
$$
\frac 1 2(Z(\hat{\beta}^{(t)})-C(\hat{\beta}^{(t)})^T\beta)^T B(\hat{\beta}^{(t)})^{-1}(Z(\hat{\beta}^{(t)})-C(\hat{\beta}^{(t)})^T\beta)+n\sum_{l=1}^p
p_{\lambda}(|\beta_l|).
$$

\item Inner loop: By the coordinate
descent method (Friedman et al, 2007) and the method of one-step local
linear approximation (Zou and Li, 2008), we obtain the update of the estimate
sequentially for coordinates $l=1,\dots,p$:
$$
\beta_l^{(t+1)}=\frac{S(u_l-\sum_{l'\neq l} (D_{ll'}\hat{\beta}^{(t+1)}_{l'},np'_{\lambda}(|\beta_l^{(t)}|))_{+}}{D_{ll}},
$$
where $S(z,\lambda)=\text{sign}(z)(|z|-\lambda)_{+}$ is the soft-thresholding operator. 
\end{itemize}
Throughout Step 1, $\alpha$ is fixed at the initial value $\hat{\alpha}_A.$ At the $t$th update of the outer loop, we calculate $\hat{\beta}^{t}$ based on the updated $Z,$ $C$ and $B.$ Nested within the outer loop, we cycle through the coordinates of $\hat{\beta}_l^{(t)},$ $l=1\dots,p,1\dots,p,\dots$ until convergence.  Based on this new $\hat{\beta}^{t},$ we update $Z,$ $C$ and $B$ and proceed to the $(t+1)$th update of the outer loop.

\item Step 2: Selection and estimation for the association model using penalized estimating equation (\ref{equation2}).

\begin{itemize}
\item Outer loop: Based on the current estimate $\hat{\alpha}^{(t)},$ we compute the updated estimate $\hat{\alpha}^{(t+1)}.$  We update the offsets term in Equation (\ref{offset}) using $\alpha^{(t)}$
and $\hat{\beta}.$ Then we perform the penalized offset logistic regression of $y_{ij}$
on $y_{ik}$ with a total of $\sum _{n_i}C_2$ observations. We iterate the update, $t=1,2, 3, \dots,$ until $||\hat{\alpha}^{(t)}-\hat{\alpha}^{(t+1)}||\leq \epsilon,$ the prespecified tolerance level.

\item Inner loop: We update the penalized estimator $\hat{ \alpha}$ coordinate-wise
for $m=1,\dots,q, 1,\dots,q,\dots$ until convergence:
$$
\alpha_m^{(t+1)}=\frac{S(
u^*_m-\sum_{m'\neq m} D^*_{mm'}\hat{\alpha}^{(t+1)}_{m'},np'_{\lambda}(|\alpha_m^{(t)}|))_{+}}{D^*_{mm}}.
$$
\end{itemize}
Throughout Step 2, $\beta$ is fixed at the final value $\hat{\beta}$ obtained from Step 1.
\end{itemize}

\begin{center}
6. {N\scriptsize UMERICAL STUDIES}
\end{center}
We conduct simulation studies to assess the performance of the
proposed methods under various circumstances. 
Binary response vectors $y_i=(y_{i1},y_{i2},\dots,y_{in_i})^T$ are
generated from the joint probability function
$$
f(y_{i1},y_{i2},\dots,y_{in_i})=\prod_{j=1}^{n_i}
\mu_{ij}^{y_{ij}}(1-\mu_{ij})^{1-y_{ij}}\{1+\Sigma_{j<k}
\rho_{ijk}
\frac{y_{ij}-\mu_{ij}}{\sqrt{v_{ij}}}\frac{y_{ik}-\mu_{ik}}{\sqrt{v_{ik}}}
\},
$$
where $\mu_{ij}=E(Y_{ij}|X_i,Z_i)$, $v_{ij}=\mu_{ij}(1-\mu_{ij}),$ $\mu_{ijk}=P(Y_{ij}=1,Y_{ik}=1|X_i, Z_i)$ and $\rho_{ijk}$ is the correlation coefficient of $Y_{ij}$ and
$Y_{ik}$ given by
$\rho_{ijk}=(\mu_{ijk}-\mu_{ij}\mu_{ik})/\sqrt{v_{ij}v_{ik}}.$ The
mean is modelled as
$$
\text{logit} \mu_{ij}=\beta_0+\sum_{l=1}^{d_1} \beta_{xl} x_{ijl}+
\sum_{l'=1}^{d_2} \beta_{zl}z_{ijl'}.
$$
The regression coefficients are set as $\beta_0=-1.6,$ $\beta_x=(3.0,0,0,1.5,0)^T,$ and
$\beta_z=(0,0,-1.5,0,0)^T.$ Covariates $x_{ij}$ are generated from a $d_1$-multivariate normal
distribution $N(\mu_x, \Sigma_x),$ and covariates $z_{ij}$ are
independently generated from a $d_2$-multivariate normal
distribution $N(\mu_z,\Sigma_z).$ We set $d_1=d_2=5,$
$\mu_x=0.5\times 1_{d_1},$  $\mu_z=-0.2\times 1_{d_2},$  $\Sigma_x$ as a $d_1\times d_1$
matrix with $(i,j)$ element being
$\sigma_{xij}=\sigma_x^2\rho_x^{|i-j|},$ $\Sigma_z$ as a
$d_2\times d_2$ matrix with $(i,j)$ element being
$\sigma_{zij}=\sigma^2_{z}\rho_z^{|i-j|},$ $\sigma_x=\sigma_z=1,$
and $\rho_x=\rho_z=0.5.$ 

The
odds ratio $\phi_{ijk}$ between $Y_{ij}$ and $Y_{ik}$ is modelled
as
$$
\log \phi_{ijk}=\alpha_0+\sum_{l=1}^{d_3}
\alpha_{wl}w_{ijkl}+\sum_{l'=1}^{d_4} \alpha_{zl'}v_{ijkl'}.
$$
The regression coefficients are set as $\alpha_0=0.693,$ $\alpha_w=(0.3,-0.3,0,0,0)^T,$
and $\alpha_v=(0,0,0,0,0)^T.$ Covariates $w_{ijkl}$ are generated from a $d_3$-multivariate normal
distribution $N(\mu_w, \Sigma_w),$ and covariates $v_{ijkl'}$ are
independently generated from a $d_4$-multivariate normal
distribution $N(\mu_v,\Sigma_v).$  We set  $d_3=d_4=5,$
$\mu_w=0.5\times 1_{d_3},$ $\mu_v=-0.2\times 1_{d_4},$
 $\Sigma_w$ as a $d_3\times d_3$
matrix with $(i,j)$ element being
$\sigma_{wij}=\sigma_w^2\rho_x^{|i-j|},$ $\Sigma_z$ as a
$d_4\times d_4$ matrix with $(i,j)$ element being
$\sigma_{vij}=\sigma^2_{v}\rho_z^{|i-j|},$ $\sigma_w=\sigma_v=1,$
and $\rho_w=\rho_v=0.5.$

The cluster size $n_i$ is set as $ 5.$ We carry out three analyses of variable selection for three different scenarios:
(1). We select the significant mean parameters with no selection conducted for
the association parameters; (2). We focus on variable selection for the association
parameters with no selection conducted for the mean parameters; (3). We conduct variable selection for both the mean and the association
parameters. In Analysis 1, we only perform the first stage penalized estimation on the mean parameters with the second stage being omitted.  To select the optimum tuning parameter, we used modified BIC information criterion, which takes the form: $
\text{BIC}=\{ \sum_i U_{i\beta}
(\hat{\beta},\hat{\alpha}_{A})\}^T\{\sum_i U_{i\beta}
(\hat{\beta},\hat{\alpha}_{A})U_{i\beta}
(\hat{\beta},\hat{\alpha}_{A})^T\}^{-1}\{ \sum_i
U_{i\beta} (\hat{\beta},\hat{\alpha}_{A})\}
+\log n \{\sum_{l=1}^{p} I(\hat{\beta}_{l}\neq
0)\}.
$     In Analysis 2, we perform the second stage penalized estimation on the association parameters with the first stage having zero penalty and hence the mean parameter estimate equal to the unpenalized ALR estimate. The BIC takes the following form:  $
\text{BIC}=\{ \sum_i U_{i\alpha}
(\hat{\beta}_A,\hat{\alpha})\}^T\{\sum_i U_{i\alpha}
(\hat{\beta}_A,\hat{\alpha})U_{i\alpha}
(\hat{\beta}_A,\hat{\alpha})^T\}^{-1}\{ \sum_i
U_{i \alpha} (\hat{\beta}_A,\hat{\alpha})\}
+\log n \{\sum_{m=1}^{q} I(\hat{\alpha}_{m}\neq
0)\}.$  In Analysis 3, the corresponding BIC takes the following form:
\begin{align}
\begin{split}
\text{BIC}=&\{ \sum_i U_{i\beta}
(\hat{\beta},\hat{\alpha}_{A})\}^T\{\sum_i U_{i\beta}
(\hat{\beta},\hat{\alpha}_{A})U_{i\beta}
(\hat{\beta},\hat{\alpha}_{A})^T\}^{-1}\{ \sum_i
U_{i\beta} (\hat{\beta},\hat{\alpha}_{A})\}
\\
&+ \{ \sum_i U_{i\alpha}
(\hat{\beta},\hat{\alpha})\}^T\{\sum_i U_{i\alpha}
(\hat{\beta},\hat{\alpha})U_{i\alpha}
(\hat{\beta},\hat{\alpha})^T\}^{-1}\{ \sum_i
U_{i \alpha} (\hat{\beta},\hat{\alpha})\}
\\
&+\log n  \{\sum_{l=1}^{p} I(\hat{\beta}_{l}\neq
0)+\sum_{m=1}^{q} I(\hat{\alpha}_{m}\neq
0)\}.
\end{split}
\end{align}
The final model is chosen with the smallest BIC value. For Analysis 1, sample sizes with $n=200,$ $500,$  and $1000$
are considered, while for Analysis 2 and 3, we choose sample sizes with $n=500,$  $1000$ and $2000$
are considered. One hundred data sets are simulated for each
parameter setting.

Table 1 summarizes the performance of the variable selection procedure for the mean parameters, i.e., the results for Analysis 1.  There are 11 mean parameters, 4 of which are nonzero and 7 of which are set to zero. The procedure is implemented with both LASSO and SCAD penalties for comparison. The optimum tuning parameter $\lambda$ is chosen through the BIC criterion.  Based on the output of the penalized estimates, average positive selections ($PS$) and average false discoveries ($FD$) are calculated. The penalized estimator using the SCAD
penalty together with the BIC information criterion is able to achieve
average $PS=4,$  the true number of nonzero coefficients. It also achieves small average $FD=0.36$ for $n=500$ and average $FD=0.02$ for $n=1000.$ In comparison, the penalized estimator using the LASSO
penalty with the BIC information criterion is also able to achieve
average $PS=4.$ However, it consistently has higher average false discoveries: average $FD=2.885$ for $n=500$ and average $FD=2.719$ for $n=1000.$

Table 2 summarizes the performance of the variable selection procedure for the association parameters, i.e., the results for Analysis 2.  There are 11 association parameters, 3 of which are nonzero and 8 of which are set to zero.  The penalized estimator using the SCAD
penalty together with the BIC information criterion is able to achieve
positive selections close to the true number of nonzero coefficients: $PS=2.4$ for $n=500$ and  $PS=2.98$ for $n=2000$ . It has average $FD=3.27$ for $n=500$ and average $FD=2.91$ for $n=2000.$ This result shows that not surprisingly, compared to the selection of mean parameters, variable selection of association parameters has higher number of false discoveries. The penalized estimator using the LASSO
penalty with the BIC information criterion is also able to achieve similar positive selections: $PS=2.53$ for $n=500$ and  $PS=2.93$ for $n=2000$ . However, compared to the SCAD penalty, the estimator with the LASSO penalty consistently has higher average false discoveries: average $FD=4.65$ for $n=500$ and average $FD=4.67$ for $n=2000.$

Table 3 summarizes the performance of the variable selection procedure for the joint selection of mean and association parameters, i.e., the results for Analysis 3.  There are 22 parameters in total, 7 of which are nonzero and 15 of which are set to zero.
The penalized estimator with the SCAD penalty maintains its satisfactory performance. Its average $PS$
is always close to the true number of nonzero coefficients: $PS=6.77$ for $n=500$ and  $PS=6.99$ for $n=2000.$  It has average $FD=7.13$ for $n=500$ and average $FD=3.39$ for $n=2000.$ Among the false discoveries, most are attributed from falsely identified association parameters.  Its FD
improves when sample sizes increases.

In summary, the penalized estimator based on the SCAD penalty tends to outperform the penalized estimator based on the LASSO penalty. The selection procedures for the mean model perform better than those for the association model. Both high sensitivity and selectivity are achieved for the selection of mean parameters, while for the selection of association parameters, only sensitivity is high but selectivity is moderate. This phenomenon is consistent with the estimation results. For example, based on one random data set with $n=2000$ observations generated from the setting above , the 8 unpenalized estimates for the zero association parameters have values range from $0.006$ to $0.27,$ and assciated standard errors ranging from $0.09$ to $0.11$; in contrast, the 7 unpenalized estimates for zero mean parameters have values range from $0.001$ to $0.04$ and associated standard errors range from $0.04$ to $0.07.$ 

It is also noted that when the modified BIC is used with the LASSO estimator, the variable selection exhibits high false discovery rates. This is because the LASSO penalty increases linearly with the size of the estimator, which leads to large biases for large nonzero parameters. This forces the procedure to select smaller tuning parameter to make the first term in BIC criterion small. By comparing the average optimum tuning parameters in Tables 1, 2 and 3, we can see the the modified BIC consistently selects smaller tuning parameters and hence has higher false discoveries for the LASSO estimator than the SCAD estimator.

\begin{center}
7. {D\scriptsize ATA ANALYSIS}
\end{center}
We apply the proposed method to analyze the data arising from a smoking cessation study ( Gruder et al., 1993, Hedeker \& Gibbons, 2006). The data
contain repeated measurements at four different time points for
489 individuals. The outcome is the dichotomous measurement representing
whether or not an individual has quit smoking. Data were collected at four telephone
interviews: post-intervention, and 6, 12, and 24 months later. The subjects are divided into four groups by the treatments they received: 
(1). randomized to the control condition; 
(2). randomized to receive a group treatment, but never showed
up to the group meetings;
(3). randomized to and received group meetings; and
(4). randomized to and received enhanced group meetings. 
In the analysis, these four groups were compared using Helmert contrasts. The three treatment contrasts were denoted as H1, H2, and H3. Additional covariates that may influence the probability of smoking abstinence include time (T),  a race indicator (race), an indicator of TV intervention (tv) and an
indicator of manual intervention (manual). 

We consider a mean model with the main covariates effects and several interaction terms included:
\begin{align}
\begin{split}
\text{logit} \mu_{ij} &=\beta_0+\beta_1T_j+\beta_2 T_j^2+\beta_3 H1_i+\beta_4 H2_i +\beta_5 H3_i+ \beta_6 race_i +\beta_7 tv_i +\beta_8 manual_i\\
&+\beta_9 T_j\times H1_i+\beta_{10} T_j\times H2_i+\beta_{11} T_j\times H3_i.
\end{split}
\end{align}
For the association
structure, we examine the empirical correlation matrix among the
four repeated measurements, and it shows that  the correlation becomes weaker
as the time interval gets larger. Therefore, for the association structure, we consider the model 
$$\log
\phi_{ijj'}=\alpha_0+\alpha_1 |T_j-T_{j'}|+\alpha_2
|T_j-T_{j'}|^2.$$

We perform penalized estimation and variable selection on
both mean and association parameters simultaneously. Using the
modified BIC, we obtain the optimum tuning parameter
$\lambda=0.09.$ Table 4 summarizes the result of the penalized estimators  obtained from the proposed two-stage penalized estimating equations. Through variable selection,
7 out of the 12 mean parameters are estimated to be
nonzero, and 2 of the 3 association parameters are estimated to be
nonzero. Both the linear and quadratic time parameters are significant, which implies the overall change in smoking abstinence involves linear and quadratic time effect. As the estimate for the linear effect is negative and the estimate for the quadratic effect is positive, a decelerating negative trend is suggested. Among the three treatment contrasts, H2 is estimated to have zero effect. This implies that whether or not showing up to the group meetings is not associated with subsequent cessation. H1 is found to be significant, which implies whether randomization to group versus to control had effect on subsequent cessation. H3 is also found to be significant, indicating that the type of meetings influenced the outcome. The race indicator is not found to be significant, indicating that there is no evidence of suggesting different effects between the white population and other ethnic groups. The TV intervention and manual instructions are found to significantly increase the probability of smoking abstinence. All of the interaction effects between time and treatments are estimated to be zero.
 For the association structure, it seems
that the log odds ratio does decrease linearly with the time
difference. The quadratic term of the time difference
is found to be insignificant in modelling the log odds ratio. Compared to the unpenalized estimates by the ALR method, the two methods seem in good agreement. The parameters that are thresholded to zero by our penalized method are those with large p-values generated by ALR method. As the tuning parameter is selected by the data driven information criterion, our method achieves the variable selection without employing any preset significance level in contrast to traditional variable selection methods.

\begin{center}
8. {C\scriptsize ONCLUSION}
\end{center}
We present a hierarchical two-stage procedure to perform simultaneous selection and estimation on
generalized estimating equations for which both mean and association structures
are modelled. The asymptotic
behavior of the penalized estimates has been established.
The numerical results and data analysis illustrate the practical
utility of the proposed method and demonstrate satisfactory performance. The proposed method can be modified to deal with data with more complex association structures. For example, clustered data frequently arise in longitudinal studies. Common examples include school-based longitudinal studies and community-based longitudinal survey. For such data, association structures typically involve three types of correlation: the clustering effect among subjects within cluster at a given time point, serial correlation of replicate measurements within subjects, and mixed effects of different subjects within the same cluster at different time points. We may adapt the discussion of Yi \& Cook (2002) to develop a simultaneous selection and estimation method to handle such data.

\begin{center}
{A\scriptsize CKNOWLEDGEMENT}
\end{center}
This research is supported by the Canadian National Science and Engineering Council grants to Gao and Yi.
 \begin{center}
 {A\scriptsize PPENDIX A}
\end{center}
Assumptions of regularity conditions:
1) $H(\beta_0,\alpha_0)$ and $V(\beta_0,\alpha_0)$ are both finite and positive definite.
2) For any $\epsilon>0,$ there exist $\eta$ and $\delta$ such that for any $|\beta-\beta_0|<\eta,$ and $|\alpha-\alpha_0|<\delta$,
we have
$$
|H(\beta_0,\alpha_0)|(1-\epsilon)\leq |H(\beta,\alpha)|\leq |H(\beta_0,\alpha_0)|(1+\epsilon),
$$
and
$$
|V(\beta_0,\alpha_0)|(1-\epsilon)\leq |V(\beta,\alpha)|\leq |V(\beta_0,\alpha_0)|(1+\epsilon),
$$
and the inequalities hold true componentwise.
3) For $n=1,$ the covariance matrices $\text{Var}_{\beta_0,\alpha_0}$ $(U_{\beta\beta}(\beta,\beta)),$ $\text{Var}_{\beta_0,\alpha_0}$ $(U_{\beta\beta}(\beta,\alpha)),$ $\text{Var}_{\beta_0,\alpha_0}$ $(U_{\beta\beta}(\alpha,\beta)),$ $\text{Var}_{\beta_0,\alpha_0}$ $(U_{\beta\beta}(\alpha,\alpha))$ are finite and positive definite.
 
\begin{center}
 {A\scriptsize PPENDIX B: PROOF OF THE THEOREMS}
\end{center}

\noindent \textbf{Proof to Theorem \ref{thm1}}
\begin{proof}
Consider a ball $||\beta-\beta_0||\leq Mn^{-\frac 1 2}$ for some finite $M.$ Applying Taylor
Expansion, we obtain:
\begin{align}
\begin{split}
G_j(\beta,\hat{\alpha}_A)&=U_{\beta_j}(\beta,\hat{\alpha}_A)-n
p'_{\lambda_n}(|\beta_j|)\text{sign}(\beta_j)\\
&=U_{\beta_j}(\beta_0,\hat{\alpha}_A)+\sum_{l=1}^p (\beta_{l}-\beta_{l0}) \partial U_{\beta_j}(\beta^*,\hat{\alpha}_A)/\partial \beta_l -n
p'_{\lambda_n}(|\beta_j|)\text{sign}(\beta_j),
\end{split}
\end{align}
for $j=1,\dots,p,$ and some $\beta^*$ between $\beta$ and $\beta_0.$
As $E(U_{\beta_j}(\beta_0,\hat{\alpha}_A))=0,$  $U_{\beta_j}(\beta_0,\hat{\alpha}_A)=O_p(n^{\frac 1 2}).$
As $|\beta^*-\beta|\leq M n^{-\frac 1 2}$ and $\hat{\alpha}_A-\alpha_0=O_p(n^{-\frac 1 2}),$ $U_{\beta \beta}(\beta^*,\hat{\alpha}_A)=O_p(n)$ componentwise.
First we consider $j\in s^c$.
Because
$
\text{lim inf}_{n\rightarrow \infty}\text{lim
inf}_{\theta\rightarrow 0+}p'_{\lambda_n}(\beta)/\lambda_n>0,
$
and $\lambda_n\rightarrow 0,$ and $\sqrt{n}\lambda_n \rightarrow
\infty$ as $n\rightarrow \infty,$ the third term dominates the the first two terms. Thus the sign of $G_j(\beta,\hat{\alpha}_A)$ is completely determined by the
sign of $\beta_j.$ This entails that inside this
$Mn^{-1/2}$ neighborhood of $\beta_0$, $G_j(\beta,\hat{\alpha}_A)>0,$ when $\beta_j<0$ and
$G_j(\beta,\hat{\alpha}_A)<0,$ when $\beta_j>0.$ Therefore for any $\beta$ inside this ball, if $\beta_{j}=0,$ it solves the equation $G_j(\beta,\hat{\alpha}_A)=0$ with probability tending to one.  This entails for any finite constant $M,$ and any sequence of $\hat{\beta},$ such that $\hat{\beta}_{j}=0,$ and $||\hat{\beta}-\beta_0||< M n^{-\frac 1 2},$ we have
$
\lim_{n\rightarrow \infty} P(G_j(\hat{\beta},\hat{\alpha}_A)=0)=1,
$ for all $j\in s^c.$

Next we consider $j\in s.$  For $n$ large enough and $\beta_{j0}\neq 0,$
$p'_{\lambda}(|\beta_{l0}|)=0$ and
$p''_{\lambda}(|\beta_{l0}|)=0.$ For notational convenience, denote $W_{\beta\beta}^*= -U_{\beta\beta}(\beta^*,\hat{\alpha}_A),$ $H=H(\beta_0,\alpha_0),$ and $H^*=H(\beta^*,\hat{\alpha}_A).$ Let $\beta_s$ denote the sub-vector of $\beta,$ $G_s$ denote the subset of penalized equations in $G,$ $U_{\beta_s}$ denote the subset of equations in $U_{\beta}$  and $W_{ss}^*$ denote the sub-matrix of $W_{\beta\beta}^*,$ $H_{ss}$ denote the sub-matrix of $H,$ and $H_{ss}^*$ denote the sub-matrix of $H^*$ corresponding to the subset of indices $s.$ This leads to the formulation:
\begin{align}
\label{GS11}
G_s(\beta,\hat{\alpha}_A)=U_{\beta_s}(\beta_0,\hat{\alpha}_A)-W_{ss}^*(\beta_s-\beta_{s0}).
\end{align}
Because of weak law of large numbers, $W_{ss}^*=n(H_{ss}+H_{ss}^*-H_{ss}+o_p(1)).$ Let $1_s$ denote a vector of ones with length equal to the cardinality of $s.$  For a given $\epsilon,$ choose $n$ sufficiently large, so that $|(H_{ss}^*-H_{ss})H_{ss}^{-1} 1_s|\leq \epsilon 1_s$ componentwise. Choose $\beta_s-\beta_{s0}=H_{ss}^{-1}1_sCn^{-\frac 1 2}$ so that $||\beta_s-\beta_{s0}||\leq Mn^{-\frac 1 2}.$ This entails $W_{ss}^*(\beta_s^n-\beta_{s0})=n(H_{ss}+H_{ss}^*-H_{ss}+o_p(1))H_{ss}^{-1}1_sCn^{-\frac 1 2}=Cn^{\frac 1 2}(1_s(1+\epsilon)+o_p(1)).$ By choosing $M$ large enough, hence $C$ large enough, we have the second term of equation (\ref{GS11}) dominating the first term. Therefore, if $\beta_s-\beta_{s0}=H_{ss}^{-1}1_sCn^{-\frac 1 2},$ $G_s(\beta,\hat{\alpha}_A)<0,$ and similarly if $\beta_s-\beta_{s0}=-H_{ss}^{-1}1_sCn^{-\frac 1 2},$ $G_s(\beta,\hat{\alpha}_A)>0.$ Because $G_s(\beta,\hat{\alpha}_A)$ is continuous on this compact set $\beta=\{(\beta_{s0}+H_{ss}^{-1}1_sdn^{-\frac 1 2},0_{p-p'})^T; -C\leq d \leq C\},$ where $0_{p-p'}$ denotes a vector of zeros of length $p-p',$ and $p'$ is the cardinality of $s.$Therefore, there exists a $\hat{\beta}$ that lies in this compact set and  $G_s(\hat{\beta},\hat{\alpha}_A)=0.$

\end{proof}

\noindent \textbf{Proof to Theorem \ref{thm2}}
\begin{proof}
Based on Taylor expansion presented in Proof to Theorem 1, we have
\begin{align}
0=G_s(\hat{\beta},\hat{\alpha}_A)=U_{\beta_s}(\beta_0,\hat{\alpha}_A)-W_{ss}^*(\hat{\beta}_s-\beta_{s0})-nb_1-n\Sigma_1^*(\hat{\beta}_s-\beta_{s0}),
\end{align}
where $\Sigma_1^*=\text{diag}\{p''_{\lambda_n}(|\beta^*_{10}|),\dots,p''_{\lambda_n}(|\beta^*_{s0}|)\},
$ and $\beta^*$ lies between $\hat{\beta}$ and $\beta_0.$
This entails
$$
(W_{ss}^*+n\Sigma_1^*)^{-1}(\hat{\beta}_s-\beta_{s0})=U_{\beta_s}(\beta_0,\hat{\alpha}_A)-nb_1.
$$
As $\hat{\beta}\rightarrow \beta_0$ in probability, $\frac 1 n W_{ss}^*\rightarrow H_{ss}$ in probability and
$\Sigma_1^*\rightarrow \Sigma_1$ in probability.
It can also be shown that
$$
\frac 1 {\sqrt{n}} U_{\beta_s}(\beta_0,\hat{\alpha}_A)= \frac 1 {\sqrt{n}} U_{\beta_s}(\beta_0,\alpha_0)+\frac 1 {\sqrt{n}} (\hat{\alpha}_A-\alpha_0)\frac {\partial U_{\beta_s}(\beta_0,\alpha^*)}{\partial \alpha}.$$ As $\sqrt{n}(\hat{\alpha}_A-\alpha_0)$ is bounded in probability and $\frac 1 n \frac {\partial U_{\beta_s}(\beta_0,\alpha^*)}{\partial \alpha}\rightarrow E\frac {\partial U_{\beta_s}(\beta_0,\alpha_0)}{\partial \alpha}=0$  in probability,
the limiting distribution of $
\frac {1} {\sqrt{n}} U_{\beta_s}(\beta_0,\hat{\alpha}_A)$ is $N\{0,V_{ss}\}.$
According to Slutsky's theorem, we
have$
\sqrt{n}(H_{ss}+\Sigma_1)\{\hat{\beta_s}-\beta_{s0}+(H_{ss}+\Sigma_1)^{-1}b_1\}\rightarrow
N\{0,V_{ss})\}.$
\end{proof}

\noindent \textbf{Proof to Theorem \ref{thm3}}
\begin{proof}

Consider a ball $||\alpha-\alpha_0||\leq Mn^{-\frac 1 2}$ for some finite $M.$ Applying Taylor
Expansion, we obtain:
\begin{align}
\label{tay1}
\begin{split}
J_j(\hat{\beta},\alpha)&=U_{\alpha_j}(\hat{\beta},\alpha)-n
p'_{\lambda_n}(|\alpha_j|)\text{sign}(\alpha_j)\\
&=U_{\alpha_j}(\beta_0,\alpha_0)+\sum_{l=1}^p (\hat{\beta}_{l}-\beta_{l0}) \partial U_{\alpha_j}(\beta^*,\alpha^*)/\partial \beta_l\\
&+\sum_{m=1}^q (\alpha_{m}-\alpha_{m0}) \partial U_{\alpha_j}(\beta^*,\alpha^*)/\partial \alpha_m -n
p'_{\lambda_n}(|\alpha_j|)\text{sign}(\alpha_j),
\end{split}
\end{align}
for some $\beta^*$ between $\hat{\beta}$ and $\beta_0$ and some $\alpha^*$ between $\alpha$ and $\alpha_0.$
 As $E(U_{\alpha_j}(\beta_0,\alpha_0))=0,$ $U_{\alpha_j}(\beta_0,\alpha_0)=O_p(n^{\frac 1 2}).$
 As $|\hat{\beta}-\beta|= O_p( n^{-\frac 1 2})$ and $\alpha-\alpha_0\leq Mn^{-\frac 1 2},$ $U_{\alpha \beta}(\beta^*,\alpha^*)=O_p(n)$ and $U_{\alpha \alpha}(\beta^*,\alpha^*)=O_p(n)$ componentwise.
First consider $j\notin v.$ Note that the first three terms are all of order $O_p(n^{1/2}).$
As
$
\text{lim inf}_{n\rightarrow \infty}\text{lim
inf}_{\alpha\rightarrow 0+}p'_{\lambda_n}(\alpha)/\lambda_n>0,
$
and $\lambda_n\rightarrow 0,$ and $n^{\frac 1 2}\lambda_n \rightarrow
\infty$ as $n\rightarrow \infty,$ the four term dominates the the first three terms. Thus the sign of $J_j(\hat{\beta},\alpha)$ is completely determined by the
sign of $\alpha_j.$ This entails that inside this
$Mn^{-1/2}$ neighborhood of $\alpha_0$, $J_j(\hat{\beta},\alpha)>0,$ when $\alpha_j<0$ and
$J_j(\hat{\beta},\alpha)<0,$ when $\alpha_j>0.$ Therefore for any $\alpha$ inside this ball, if $\alpha_{j}=0,$ it solves the equation $J_j(\hat{\beta},\alpha)=0$ with probability tending to one.  This entails for any finite constant $M,$ and any sequence of $\hat{\alpha},$ if $\hat{\alpha}_{j}=0,$ and $||\hat{\alpha}-\alpha_0||< M n^{-\frac 1 2},$ we have
$
\lim_{n\rightarrow \infty} P(J_j(\hat{\beta},\hat{\alpha})=0)=1,
$ for all $j\notin v.$

Next we consider $j\in v.$ It is known that for $n$ large enough, and $\alpha_{j0}\neq 0,$
$p'_{\lambda}(|\alpha_{j0}|)=0$ and
$p''_{\lambda}(|\alpha_{j0}|)=0.$   Let $W^*_{\alpha\alpha}=- U_{\alpha\alpha}(\beta^*,\alpha^*),$ and $W^*_{vv}$ denote the sub-matrix of $W^*_{\alpha\alpha}.$ Let $W^*_{\alpha\beta}=-U_{\alpha\beta}(\beta^*,\alpha^*),$ and $W^*_{vs}$ denote the sub-matrix of $W^*_{\alpha\beta}$ corresponding to indices set $v$ and $s.$ Then Equation (\ref{tay1}) can be expressed as
\begin{align}
\label{GS}
J_v(\hat{\beta},\alpha)=U_{\alpha_v}(\beta_0,\alpha_0)-W_{vs}^{*}(\hat{\beta}_{s}-\beta_{s0})-W_{vv}^{*}(\alpha_{v}-\alpha_{v0}).
\end{align}
From the proof to Theorem 2, we have $(\hat{\beta}_{s}-\beta_{s0})=H_{ss}^{-1}U_{\beta_s}(\beta_0,\alpha_0)+o_p(1).$ Furthermore, we have
$W_{vs}^{*}\rightarrow H_{vs}$ in probability and $W_{vv}^{*}\rightarrow H_{vv}$ in probability.
Therefore, Equation (\ref{GS}) can be simplified as
\begin{align}
\label{GS1}
J_v(\hat{\beta},\alpha)=U_{\alpha_v}(\beta_0,\alpha_0)-H_{vs}H_{ss}^{-1}U_{\beta_s}(\beta_0,\alpha_0)-H_{vv}(\alpha_{v}-\alpha_{v0})+o_p(1).
\end{align}
We choose $\alpha_{v}-\alpha_{v0}=H^{vs}C_1n^{\frac 1 2}+H^{vv}C_2n^{\frac 1 2}.$ Because $H_{vs}H_{ss}^{-1}+H_{vv}H^{vs}=0,$
the left side of the equation can be expressed as
$$
-H_{vs}H_{ss}^{-1}(U_{\beta_s}(\beta_0,\alpha_0)-C_1 n^{\frac 1 2})+(U_{\alpha_v}(\beta_0,\alpha_0)-C_2n^{\frac 1 2}).
$$
Because both $U_{\beta_s}(\beta_0,\alpha_0)$ and $U_{\alpha_v}(\beta_0,\alpha_0)$ are of order $O_p(n^{\frac 1 2}),$ we can choose  $M$ large enough, hence $C_1$ and $C_2$ large enough so that both  $U_{\beta_s}(\beta_0,\alpha_0)$ and $U_{\alpha_v}(\beta_0,\alpha_0)$ are dominated by $C_1 n^{\frac 1 2},$ $C_2 n^{\frac 1 2}.$  Furthermore, $J_v(\hat{\beta},\alpha)$ is continuous on this compact set $\{(H^{vs}\rho C_1n^{\frac 1 2}+H^{vv}\rho C_2n^{\frac 1 2},0_{q-q'})^T;-1\leq \rho \leq 1\},$ where $0_{q-q'}$ denotes a vector of zeros of length $q-q',$ and $q'$ is the cardinality of $v.$ The sign of $J_v(\hat{\beta},\alpha)$ is opposite when $\rho=-1,$ and $1.$Therefore, there exist a $\hat{\alpha}$ that lies in this compact set and  $J_v(\hat{\beta},\hat{\alpha})=0.$

\end{proof}

\noindent \textbf{Proof to Theorem \ref{thm4}}
\begin{proof}
Based on the proofs to Theorem 1, 2 and 3, Taylor expansions of the penalized estimating equations at the penalized estimators can be expressed as:
$$
(\hat{\beta_s}-\beta_{s0})(W_{ss}^*+n\Sigma_1)+(\hat{\alpha}_A-\alpha_0)W_{\beta\alpha}^*=U_{\beta_s}(\beta_0,\alpha_0)-nb_1,
$$
and
$$
(\hat{\beta_s}-\beta_{s0})(W_{vs}^*)+(\hat{\alpha}_s-\alpha_{s0})(W_{vv}^*+n\Sigma_2)=U_{\alpha_v}(\beta_0,\alpha_0)-nb_2.
$$
Because $\frac 1 n W^*_{ss}\rightarrow H_{ss}$ in probability, $\frac 1 n W^*_{\beta\alpha}\rightarrow H_{\beta\alpha}=0$ in probability, $\frac 1 n W^*_{vv}\rightarrow H_{vv}$ in probability, we have
$$
(\hat{\beta}_s-\beta_{s0})(H_{ss}+\Sigma_1)=U_{\beta_s}(\beta_0,\alpha_0)/n-b_1+o_p(1),
$$
and
$$
(\hat{\beta}_s-\beta_{s0})(H_{vs})+(\hat{\alpha}_v-\alpha_{v0})(H_{vv}+\Sigma_2)=U_{\alpha_v}(\beta_0,\alpha_0)/n-b_2+o_p(1).
$$
This implies
$$
\sqrt{n}\left(
          \begin{array}{c}
            \hat{\beta}_s-\beta_0 \\
            \hat{\alpha}_v-\alpha_0 \\
          \end{array}
        \right)\left(
                                      \begin{array}{cc}
                                        H_{ss}+\Sigma_1 & 0 \\
                                        H_{vs} & H_{vv}+\Sigma_2 \\
                                      \end{array}
                                    \right)=\sqrt{n}\left(
          \begin{array}{c}
            U_{\beta_s}(\beta_0,\alpha_0)/n-b_1 \\
            U_{\alpha_v}(\beta_0,\alpha_0)/n-b_2 \\
          \end{array}
        \right)+o_p(1).
$$
According to Slutsky's theorem and central limit theorem, the joint distribution of $\hat{\beta}_s$ and $\hat{\alpha}_v$ converges in distribution to the joint multivariate normal distribution.
\end{proof}

\begin{center}
 {A \scriptsize PPENDIX C: VARIANCE ESTIMATE}
\end{center}

Here we outline how to evaluate the standard errors of the
penalized estimators. The consistent estimate for the negative Hessian matrix $H$ is denoted as
$$
\hat{H}=\left(%
\begin{array}{cc}
  \sum_i C_{i}(\hat{\beta},\hat{\alpha})^T B(\hat{\beta},\hat{\alpha})^{-1}_i C_{i}(\hat{\beta},\hat{\alpha}) & 0 \\
  \sum_i T_{i}(\hat{\beta},\hat{\alpha})^T S_i(\hat{\beta},\hat{\alpha})^{-1} F_{i}(\hat{\beta},\hat{\alpha}) & \sum_i T_{i}(\hat{\beta},\hat{\alpha})^T S_i(\hat{\beta},\hat{\alpha})^{-1} T_{i}(\hat{\beta},\hat{\alpha}) \\
\end{array}%
\right),
$$
where $F_{i}=\partial \zeta_i/\partial \beta.$ The consistent
estimator of the variance matrix $V$ of the score vectors is denoted as
$$
\hat{V}=\left(%
\begin{array}{cc}
  \sum_i
U_{i\beta}(\hat{\beta},\hat{\alpha})U_{i\beta}(\hat{\beta},\hat{\alpha})^T &   \sum_iU_{i\beta}(\hat{\beta},\hat{\alpha})U_{i\alpha}(\hat{\beta},\hat{\alpha})^T \\
 \sum_i
U_{i\alpha}(\hat{\beta},\hat{\alpha})U_{i\beta}(\hat{\beta},\hat{\alpha})^T &   \sum_iU_{i\alpha}(\hat{\alpha},\hat{\alpha})U_{i\alpha}(\hat{\beta},\hat{\alpha})^T \\
\end{array}%
\right)
,
$$
with $U_{i\beta}(\hat{\beta},\hat{\alpha})=C_{i}(\hat{\beta},\hat{\alpha})^T B(\hat{\beta},\hat{\alpha})^{-1}_i A_{i}(\hat{\beta},\hat{\alpha}),$ and $U_{i\alpha}(\hat{\beta},\hat{\alpha})=T_{i}(\hat{\beta},\hat{\alpha})^T S(\hat{\beta},\hat{\alpha})^{-1}_i R_{i}(\hat{\beta},\hat{\alpha}).$
We also define
$$
\hat{\Sigma_{1}}=\text{diag}\{p_{\lambda}'(|\hat{\beta}_{1}|)/|\hat{\beta}_{1}|,\dots,p_{\lambda}'(|\hat{\beta}_{p'}|)/|\hat{\beta}_{p'}|\}
,$$ and
$$
\hat{\Sigma_{2}}=\text{diag}\{p_{\lambda}'(|\hat{\alpha}_{1}|)/|\hat{\alpha}_{1}|,\dots,p_{\lambda}'(|\hat{\alpha}_{p'}|)/|\hat{\alpha}_{q'}|\}.
$$
Then estimated covariance matrix for $\sqrt{n}\left(
          \begin{array}{c}
            \hat{\beta}_s-\beta_{s0} \\
            \hat{\alpha}_v-\alpha_{v0} \\
          \end{array}\right)$ is
$
                              \hat{B}\left(
                                      \begin{array}{cc}
                                        \hat{V}_{ss} & \hat{V}_{sv} \\
                                        \hat{V}_{vs} & \hat{V}_{vv} \\
                                      \end{array}
                                    \right)\hat{B}^T,
$ where
$\hat{B}= \left(
                                      \begin{array}{cc}
                                        \hat{H}_{ss}+\hat{\Sigma}_1 & 0 \\
                                       \hat{H}_{vs} & \hat{H}_{vv}+\hat{\Sigma}_2 \\
                                      \end{array}
                                    \right)^{-1}.$
\vskip 0.5cm

\newpage
\noindent{\large\bf References}
\begin{description}

\item Bondell, H. D., Krishna, A., \& Ghosh, S. K. (2010). Joint variable selection
for fixed and random effects in linear mixed-effects models. \emph{Biometrics} \textbf{66}, 1069-1077.

\item Cai, J., Fan, J., Li, R., \& Zhou, H. (2005). Variable selection for
multivariate failure time data. \emph{Biometrika} \textbf{92}, 303–316.

\item Carey, V., Zeger, S. \& Diggle, P. (1993). Modelling multivariate binary data with alternating logistic regressions. \emph{Biometrika} \textbf{80}, 517-526.

\item Fan, J. \& Li, R. (2001). Variable selection via nonconcave
penalized likelihood and its oracle properties. \emph{Journal of
the American  Statistical Association} \textbf{96}, 1348-60.

\item Fan, J. \& Li, R. (2002). Variable selection for Cox’s proportional
hazards model and frailty model. \emph{Annals of Statistics} \textbf{30}, 74–
99.

\item Fan, J. \& Li, R. (2004). New estimation and model selection procedures
for semiparametric modeling in longitudinal data analysis.
\emph{Journal of the American Statistical Association} \textbf{99}, 710–
723.

\item Fitzmaurice, G. M. \& Lipsitz, S. R. (1995). A Model for Binary Time
Series Data With Serial Odds Ratio Patterns. \emph{Applied Statistics} \textbf{44},
51–61.

\item Friedman, J., Hastie, T. \&  Tibshirani, R. (2008). Sparse
inverse covariance estimation with the graphical lasso.
\emph{Biostatistics} \textbf{9}, 432-441.

\item Garcia, R. I., Ibrahim, J. G., \& Zhu, H. (2010). Variable selection
for regression models with missing data. \emph{Statistica Sinica} \emph{20}, 149–
165.

\item Gruder, C.L.,  Mermelstein, R.J.,  Kirkendol, S.,  Hedeker, D.,
 Wong, S.C.,  Schreckengost, J. , Warnecke,  R.B.,  Burzette, R., and
Miller, T.Q. (1993) Effects of social support and relapse prevention training
as adjuncts to a televised smoking cessation intervention.
\emph{Journal of Consulting and Clinical psychology} \textbf{61}, 113-120.

\item He, H. \&  Yi, G. Y. (2011).
A Pairwise Likelihood Method for Correlated Binary Data
with/without Missing Observations under Generalized Partially
Linear Single-Index Models.  {\em Statistica Sinica}, \textbf{21}, 207-229.

\item Hedeker D. \&
Gibbons, R.D. (2006)  \emph{Longitudinal Data analysis}, Wiley.

\item Ibrahim, J. G., Zhu, H., Garcia, R. I. \& Guo, R. (2011). Fixed and Random Effects Selection in Mixed Effects Models. \emph{Biometrics} \textbf{67}, 495–503.

\item Lang, J. \& Agresti, A. A. (1994). Simultaneously Modeling Joint and
Marginal Distributions of Multivariate Categorical Responses. \emph{Journal of
the American Statistical Association}, \textbf{89}, 625–632.

\item Liang, K.-Y. \& Zeger, S.L. (1986). Longitudinal data analysis using generalized linear models. \emph{Biometrika} \textbf{73}, 13-22.

\item Liang, K.-Y., Zeger, S. L. \& Qaqish, B. (1992). Multivariate Regression
Analyses for Categorical Data (with discussion). \emph{Journal of the Royal
Statistical Society, Ser. B}, \textbf{54}, 3–40.

\item Lipsitz, S. R., Laird, N. M. \& Harrington, D. P. (1991). Generalized Estimating
Equations for Correlated Binary Data: Using the Odds Ratio as a
Measure of Association. \emph{Biometrika}, \textbf{78}, 153–160.

\item Meinshausen, N. \& Buhlmann, P. (2006). High-dimensional graphs and variable
selection with the lasso. \emph{Annals of Statistics} \textbf{34}, 1436-1462.

\item Molenberghs, G. \& Lesaffre, E. (1994). Marginal Modelling of Correlated
Ordinal Data Using an n-Way Plackett Distribution. \emph{Journal of the American
Statistical Association} \textbf{89}, 633–644.

\item Prentice, R. L. (1988) Correlatedbinary rgression with covariates specfic to
 each binary observation. {\em Biometrics} \textbf{44}, 1033-1048.

\item Qu, A. \& Li, R. (2006). Quadratic inference functions for varyingcoefficient
models with longitudinal data. \emph{Biometrics} \textbf{62}, 379–
391.

\item Tibshirani, R. J. (1996) Regression shrinkage and selection
via the lasso. \emph{Journal of the Royal Statistical Society,
Series B} \textbf{58}, 267-288.

\item  Yi, G. Y. \&  Cook, R. J. (2002). Marginal methods for incomplete longitudinal data arising in clusters.
{\em Journal of the American Statistical Association} \textbf{97}, 1071-1080.

\item Yi, G. Y., He, W. \& Liang, H. (2009).
Analysis of Correlated Binary Data under
Partially Linear Single-Index Logistic Models. {\em Journal
of Multivariate Analysis} \textbf{100}, 278-290.

\item Yi, G. Y., He, W. \& Liang, H. (2011).
Semiparametric Marginal and Association Regression Methods for
Clustered Binary Data.
{\em Annals of the Institute of Statistical Mathematics} \textbf{63}, 511-533.

\item Yuan, M. \& Lin, Y. (2007) Model selection and estimation in the gaussian
graphical model. \emph{Biometrika} \textbf{94}, 19-35.

\item Zou, H. \& Li, R. (2008) One-step sparse estimates in
nonconcave penalized likelihood models (with discussion).
\emph{Annals of Statistics}  \textbf{36}, 1509-1533.

\item Wang, L. (2011). GEE analysis of clustered binary data with diverging number of covariates. \emph{Annals of Statistics}, \textbf{39}, 389-417.

\end{description}

\begin{table}

 \caption{ Positive selections (PS) and false discoveries (FD) for variable selection of mean parameters with 4 nonzero coefficients and 7 zero coefficients }
\label{tab1}

\begin{center}
{\scriptsize
\begin{tabular}{ccccccc}
\hline\hline   &  &LASSO & &  &SCAD &\\\hline
 n & $\overline{\lambda}$ &PS &FD & $\overline{\lambda}$ &PS &FD\\
\hline
200 & 0.019   & 4   & 3.690   &    0.064  &       4  & 1.620 \\
    & (0.011) & (0) & (1.830) &   (0.012) &      (0) &(1.135) \\

500 & 0.015   & 4   & 2.885   &    0.058  &       4  & 0.360 \\
    & (0.006) & (0) & (1.486) &   (0.015) &      (0) &(0.578) \\

1000 & 0.011   & 4   & 2.719   &    0.053  &       4  & 0.020 \\
    & (0.005) & (0) & (1.412) &   (0.013) &      (0) &(0.140) \\


\hline\hline
\end{tabular}%
}
\end{center}
{\scriptsize
 (PS denotes the number of correctly identified nonzero coefficients;FD denotes the number of zero coefficients incorrectly estimated to be nonzero;numbers without parenthesis are average values; numbers with parenthesis are standard deviations;$\overline{\lambda}$ denotes the average optimum tuning parameter.)}
\end{table}

\begin{table}

 \caption{ Positive selections (PS) and false discoveries (FD) for variable selection of association parameters with 3 nonzero coefficients and 8 zero coefficients }
\label{tab1}

\begin{center}
{\scriptsize
\begin{tabular}{ccccccc}
\hline\hline   &  &LASSO & &  &SCAD &\\\hline
 n & $\overline{\lambda}$ &PS &FD & $\overline{\lambda}$ &PS &FD\\
\hline
500 & 0.027   & 2.530   & 4.650   &    0.074  &       2.400  & 3.270 \\
    & (0.020) & (0.688) & (1.977) &   (0.032) &      (0.752) &(2.004) \\

1000 & 0.022   & 2.730   & 4.480   &    0.054  &       2.720  & 3.120 \\
    & (0.015) & (0.510) & (2.254) &   (0.021) &      (0.570) &(2.076) \\

2000 & 0.014   & 2.930   & 4.670   &    0.040  &       2.980  & 2.910 \\
    & (0.009) & (0.256) & (2.040) &   (0.012) &      (0.141) &(1.730) \\

\hline\hline
\end{tabular}%
}
\end{center}
{\scriptsize
( PS denotes the number of correctly identified nonzero coefficients;FD denotes the number of zero coefficients incorrectly estimated to be nonzero;numbers without parenthesis are average values; numbers with parenthesis are standard deviations;$\overline{\lambda}$ denotes the average optimum tuning parameter.)}
\end{table}

\begin{table}

 \caption{ Positive selections (PS) and false discoveries (FD) for variable selection of both mean and association parameters with 7 nonzero coefficients and 15 zero coefficients }
\label{tab1}

\begin{center}
{\scriptsize
\begin{tabular}{ccccccc}
\hline\hline   &  &LASSO & &  &SCAD &\\\hline
 n & $\overline{\lambda}$ &PS &FD & $\overline{\lambda}$ &PS &FD\\
\hline
500 & 0.015   & 6.840   & 9.430   &    0.042  &       6.770  & 7.130 \\
    & (0.007) & (0.368) & (2.508) &   (0.005) &      (0.423) &(1.662) \\

1000 & 0.012   & 6.940   & 9.380   &    0.041  &       6.880  & 5.430 \\
    & (0.005) & (0.278) & (2.420) &   (0.006) &      (0.356) &(1.486) \\

2000 & 0.009   & 7.000   & 9.380   &    0.040  &       6.990  & 3.390 \\
    & (0.004) & (0) & (2.490) &   (0.006) &      (0.100) &(1.421) \\

\hline\hline
\end{tabular}%
}
\end{center}
{\scriptsize
( PS denotes the number of correctly identified nonzero coefficients;FD denotes the number of zero coefficients incorrectly estimated to be nonzero;numbers without parenthesis are average values; numbers with parenthesis are standard deviations;$\overline{\lambda}$ denotes the average optimum tuning parameter.)}
\end{table}

\begin{table}

 \caption{  Penalized estimation of smoking cessation study data set }
\label{tab1}

\begin{center}
{\scriptsize
\begin{tabular}{ccccc}
\hline\hline   &   ALR  &   & PGEE2 &
\\\hline
 variable & estimate & SE   &estimate & SE  \\
\hline
 intercept&   -1.280 &  0.140    &   1.229 & 0.123         \\
      time&   -0.850 &  0.145    &   -0.778 & 0.142     \\
  time$^2$&   0.238  &  0.044    &   0.210 & 0.043   \\
  hermert1&   0.563  &  0.211   &    0.321 & 0.175   \\
  hermert2&   0.225  &  0.149    &  0    &  0        \\
  hermert3&   0.324  &  0.140    &  0.255 & 0.121      \\
     racew&   0.295  &  0.210    &   0   &   0       \\
        tv&   0.512  &  0.201    &   0.545 & 0.199     \\
    manual&   0.516  &  0.203    &   0.505 & 0.192     \\
  timeXh1&   -0.152  &  0.077    &   0   &     0      \\
  timeXh2&   -0.080&  0.060     &   0   &      0    \\
  timeXh3&   -0.049&  0.066     &   0   &      0   \\
 \hline
 intercept&   3.582  & 0.442 &   3.103& 0.163     \\
  timediff&   -1.068  & 0.525&     -0.477& 0.078    \\
  timediff$^2$&  0.155 &  0.141  &   0   &  0        \\
\hline\hline
\end{tabular}%
}
\end{center}
{\scriptsize
(The first part of the table includes mean parameters and the second part of the table includes the association parameters.)}
\end{table}
\end{document}